\documentclass[12pt]{article}
\usepackage{amsmath,amsfonts,amssymb,amsthm,amscd}
\title{Strongly minimal pseudofinite structures}
\date{\today}
\author{Anand Pillay\thanks{Supported by  NSF grant DMS-1360702 }\\University of Notre Dame}
\newtheorem{Theorem}{Theorem}[section]

\newtheorem{Definition}[Theorem]{Definition}

\newtheorem{Corollary}[Theorem]{Corollary}
\newtheorem{Fact}[Theorem]{Fact}

\newtheorem{Problem}[Theorem]{Problem}
\newcommand{\R}{\mathbb R}   
  
\newcommand{\Z}{\mathbb Z}  
\newcommand{\N}{\mathbb N}

\begin{document}
\maketitle

\begin{abstract} 
We observe that the nonstandard finite cardinality of a definable set $X$ in a strongly minimal pseudofinite structure $D$  is a polynomial over $\Z$  in the nonstandard finite cardinality of $D$. We conclude that $D$  is {\em unimodular}  in the sense of \cite{Hrushovski}, hence also {\em  locally modular}. 
We also deduce  a {\em regularity lemma}  for graphs definable in strongly minimal pseudofinite structures, although local modularity severely restricts the examples. 
The paper is elementary, and the only surprising thing about it is that the results were not explicitly noted before.  

\end{abstract}

\section{Introduction} 

\subsection{Statement of results}
A first order theory $T$ is said to be pseudofinite if every sentence in $T$ has a finite model, and an $L$-structure $M$ is said to be pseudofinite if $Th(M)$ is.  The model-theoretic study of pseudofinite theories/structures is the same as the model theoretic study of classes of finite structures focusing on suitable  uniformities and asymptotics.  One might have thought that there is not much to say about this model theory of pseudofinite theories and structures, other than in very special cases or under strong additional assumptions, such as  smoothly approximable structures, pseudofinite fields, or the context of \cite{MS-measurable}.   In fact the convential wisdom was that in so far as model theory was relevant to (classes of) finite structures,  the logic would have to go beyond first order logic, as in so-called ``finite model theory".  However the paper \cite{Hrushovski-approximate} showed that  first order model theory, even in its tame variety, is meaningful in the  study of rather arbitrary families of finite structures, and since then  additional interest in first order pseudofinite model theory has developed.  This is the context in which we take a look at strongly minimal pseudofinite theories.  Strongly minimal theories  are the ``nicest" stable theories in various senses, and are defined/ characterized by any definable subset of the universe of a model of $T$ being finite or cofinite. As it turns out the behaviour of strongly minimal pseudofinite structures is like in pseudofinite fields but much better.  We prove:

\begin{Theorem} Let $D$ be a (saturated) pseudofinite strongly minimal structure, and let $q\in \N^{*}$ be the pseudofinite cardinality of $D$ (written  $q = |D|$). Then 
\newline 
(i) for any definable (with parameters) set $X\subseteq D^{n}$, there is polynomial $P_{X}(x)$ with integer coefficients and positive leading coefficient  such that $|X| = P_{X}(q)$. Moreover $RM(X) = degree(P_{X})$. 
\newline
(ii) In fact,  for any $L$-formula $\phi(\bar x,\bar y)$, there are a finite number $P_{1},..,P_{k}$ of polynomials over $\Z$, and formulas $\psi_{1}(\bar y),..,\psi_{k}(\bar y)$, such that
the $\psi_{i}(\bar y)$ partition $\bar y$-space, and moreover for any ${\bar b}$, $|\phi({\bar x,\bar b})(D)| = P_{i}(q)$ iff $\models \psi_{i}(\bar b)$. 
\end{Theorem}

Let us remark that part (ii) of the above is not so important for applications we mention, but nevertheless is rather striking and gives a nice ``generalization" of definability of Morley rank, or definability of (Morley rank, measure) in suitable theories. 

Recall that a strongly minimal structure is said to be unimodular if whenever $\bar a = (a_{1},..,a_{n})$ and $\bar b = (b_{1},..,b_{n})$ are each $n$-tuples of algebraically independent elements, and $\bar a$ is interalgebraic with $\bar b$, then $mult(\bar a/\bar b) = mult(\bar a/\bar b)$.  Here $mult({\bar a}/{\bar b})$ denotes the number of realizations of $tp({\bar a}/{\bar b})$ which is finite, by hypothesis.  The notion is due to Hrushovski \cite{Hrushovski}, but a clarification of the definition appears in \cite{KP}. 
The  important result in  \cite{Hrushovski} (also with an exposition in \cite{Pillay-book}) states that unimodular strongly minimal sets are locally modular (generalizing Zilber's result that $\omega$-categorical strongly minimal sets are locally modular).   Macpherson and Steinhorn \cite{MS-measurable} formulated a notion of ``measurability" for first order theories, which in \cite{KP} we called $MS$-measurability, modelled after the example of pseudofinite fields. Roughly speaking any definable set is assigned both a dimension and measure, satisfying appropriate axioms. In \cite{KP} we observed that for strongly minimal sets, unimodularity is equivalent to $MS$-measurability. Now it would be not too hard to see that the pair
$(deg(P_{X}), \ell lc(P_{X}))$  (where $\ell c$ is ``leading coefficient") witnesses $MS$-measurability of $D$ when $D$ is pseudofinite.  But in fact we will prove  unimodularity directly:
\begin{Corollary} A strongly minimal pseudofinite structure $D$ is unimodular, hence by \cite{Hrushovski} locally modular.
\end{Corollary}

Now local modularity is a severe restriction on a strongly minimal set: it implies that $D$ is either ``essentially" a vector space over a division ring $R$ (in the language of $R$-modules), or algebraic  closure inside $D$ is ``trivial" in the sense that $acl(A)= \cup_{a\in A}acl(a)$.  Maybe it is worth classifying in a suitable sense the pseudofinite possibilities.  But  any case, the regularity lemma that we deduce below, should be considered to have  limited interest in view of the very special nature of the examples.

 Szemeredi's regularity lemma is roughly speaking a combinatorial result about decomposing finite graphs $G$  into a small number of ``subgraphs" $G_{i}$ such that all but a certain ``small" exceptional set of these subgraphs $G_{i}$ have  the regularity property that the density of arbitrary subgraphs of $G_{i}$ is not too different from the density of $G$, all expressed approximately and/or asymptotically.  The regularity lemma has a nonstandard formulation for pseudofinite graphs, and also a  direct proof of such a nonstandard formulation. See Tao's blog \cite{Tao-blog} for example.   Stronger versions of the regularity lemma have been given for certain restricted classes of finite graphs. For example Tao's {\em Algebraic regularity lemma} in \cite{Tao}  concerns graphs uniformly definable in finite fields (where definability refers to the field language).  Among the improvements over the usual regularity lemma is that there are no ``exceptional pairs" (or in the notation above, exceptional subgraphs).  And again this can be and in fact {\em is}  reformulated as a theorem about definable graphs in pseudofinite fields.  A related but somewhat different strong regularity lemma appears in \cite{MS}; where one considers the class of all finite bi-partitite graphs $(V,W,E)$ such that  for a given $k$, $E$ does not have the $k$-order property (in the sense of stability).  

Now, given our (saturated) strongly minimal pseudofinite structure $D$ in language $L$ say, we have at hand not only sets definable in $D$ in the language $L$, but also ``internal" subsets of $D, D\times D,...$  in the sense of nonstandard analysis, each with their own nonstandard finite cardinality.   The following will be immediate from Theorem 1.1. 

\begin{Corollary} Let $D$ be a saturated strongly minimal  pseudofinite structure, with $|D| = q$.  Let $(V,W,E)$ be a graph definable in $D$ such that $RM(V) = d_{1}$, $RM(W) = d_{2}$  (so  $P_{V\times W}(q)$  has degree $d_{1} + d_{2}$)  and both $V, W$ have Morley degree $1$.   Then there is a polynomial $R$ with integer coefficients and of degree $< d_{1} + d_{2}$ such that either:
\newline
(i) $|E| = |V||W| - R(q)$, in which case  also for any internal subsets $A$ of $V$ and $B$ of $W$, $|E\cap (A\times B)| \geq  |A||B| -  R(q)$, or
\newline
(ii) $|E| = R(q)$, in which case for any internal subsets $A$ of $V$ and $B$ of $W$, $|E\cap (A\times B)| \leq R(q)$. 
\end{Corollary}

Now for an arbitrary definable graph $(V,W,E)$ in $D$ with $RM(V) = d_{1}$ and $RM(W) = d_{2}$ we can partition $V$ into finitely many sets $V_{1}\cup ...\cup V_{k_{1}}$ of Morley rank $d_{1}$ and Morley degree $1$  (where $d_{1}$ is the Morley degree of $V_{1}$) and likewise for $W$.   Then Corollary 1.3  holds for each of the definable graphs 
$(V_{i},W_{j}, E|(V_{i}\times W_{j}))$, and this is the strongly minimal pseudofinite  regularity lemma. 

It is routine to reformulate this regularity lemma in terms of finite graphs, namely the class of graphs uniformly definable in finite models of $Th(D)$. We leave details to the interested reader. 
% For example suppose $T = Th(D)$ is countable and let $T = \{\sigma_{i}:i<\omega}$ be an enumeration. Let $K_{n}$ be the class of finite models of $\sigma_{1}\wedge ...\wedge %\sigma_{n}$, and there is no harm arranging this such that the minimum cardinality of a structure in $K_{n}$ is $f(n)$ where $f$ is a strictly increasing function.  Then .... (is it worth it...)

The material here (at least Corollary 1.2) dates back to discussions  with various  people at the MSRI in Spring 2014, in particular Martin Bays and Pierre Simon, as well as Dugald Macpherson and Charles Steinhorn,  on the status of  the statement  ``any strongly minimal pseudofinite theory is locally modular".  The conclusion was that it does or at least {\em should}  follow from Proposition 3.1 (iv) of \cite{KP} (together with the fact that unimodularity implies local modularity). In the autumn of 2014 Alex Kruckmann asked Macpherson and myself again about the status of this problem and in working out the response I realized that in fact we have the nice polynomial counting result of  Theorem 1.1.  After noting also the (easy) regularity lemma mentioned above, it seemed worthwhile to write down the details. Anyway thanks to the above-mentioned people for the discussions and communications, as well as Sergei Starchenko and the participants in my course on pseudofinite model theory at Notre Dame, Autumn 2014.

After a first version of this preprint was circulated, Dugald Macpherson pointed out that in fact $MS$-measurability (and thus by \cite{KP},  unimodularity) of peudofinite strongly minimal sets can be deduced from Lemma 2.5 of \cite{MS-measurable}.  Also Terence Tao suggested to include some remarks on the limited nature of the examples (which we have done above).

\subsection{Preliminaries}
 Strongly minimal sets control uncountably categorical theories and more generally  most $\omega$-stable theories of finite Morley rank. The model theory of strongly minimal sets is very well-known and is the most accessible incarnation or special case of general stability theory. The reader is referred to Chapters 5 and 6 of \cite{TZ}.    But I will give below a brief  summary of what is needed in the current paper, in case there is an interested  reader from outside the subject. The basic examples of strongly minimal structures are algebraically closed fields (in the ring language), infinite vector spaces over division rings $R$ in the language of $R$-modules, and infinite free $G$-sets, in lthe language with unary function symbols for elements of the group $G$.  

I will assume for the rest of this section as well as in the proofs of the main results just basics of model theory including  notions of type, saturation, algebraic closure.  A complete $1$-sorted  theory $T$ in language $L$ is said to be {\em strongly minimal} if every definable (with parameters) subset $X$ of any model $M$ of $T$ is finite or cofinite. The condition is on definable subsets of ``$1$-space",  but conclusions are drawn about the structure of definable subsets of higher ambient spaces. 

Until otherwise mentioned $T$ denotes a complete strongly minimal theory in language $L$  and $D$ is a saturated  model. $D^{n}$ is of course just the collection of $n$-tuples from $D$, but we also allow $n=0$ in which $D^{0}$ is considered as an auxiliary point.

\begin{Definition} (i) Let $X\subseteq D^{n}$ be definable. Then $dim(X)$ is the least $k\leq n$ such that we can write $X$ as a finite union of definable sets $X_{1}\cup ..\cup X_{r}$ such that for each $i$ there is a projection $\pi_{i}$ from $D^{n}$ onto some $D^{k}$ such that the restriction $\pi_{i}|X_{i}$ of $p_{i}$ to $X_{i}$ is finite-to-one. 
\newline
(ii) Let $X\subseteq D^{n}$ be definable and of dimension $k$. Then $mlt(X)$ is the greatest natural number $m$ (if one exists) such that $X$ can be written as a disjoint union  of definable $X_{i},..,X_{m}$ such that $dim(X_{i}) = k$ for each $i$. 
\newline
(iii) By a $k$-cell we mean a definable set $X\subseteq D^{n}$ for some $n\geq k$ such that for some finite nonzero $r$ there is a projection $\pi$ from $D^{n}$ to some $D^{k}$, such that $dim(\pi(X)) = k$ and $\pi|X$ is $r$-to-$1$. 
\end{Definition} 

Clearly $dim(X)$ as defined in (i) exists. Because if $X\subseteq D^{n}$ then already the projection from $D_{n}$ to itself is one-to-one.

\begin{Fact} (i) $mlt(X)$ exists for any definable $X$.
\newline
(ii) For any $n\geq 0$, $D^{n}$ has dimension $n$ and multiplicity $1$.
\newline
(iii) Any $k$-cell  has dimension $k$. Moreover any definable $X$ is a finite disjoint of cells, i.e. of $k$-cells for possibly varying $k$.
\newline
(iv) For $X$ definable $dim(X) = 0$ iff $X$ is finite.
\end{Fact}

Note that Fact 1.5 (ii) in the case where $n = 1$ is just the definition of strong minimality. The main point is to deduce ``good behaviour" for definable sets in higher ambient spaces, and this is a good exercise for a beginner, using some basic lemmas  about algebraic closure in strongly minimal sets.  For example,  for $b_{1},..,b_{n}\in D$ and $A$ a small subset of $D$, we call $\{b_{1},..,b_{n}\}$ algebraically independent over $A$ if $b_{i}\notin acl(A,b_{1},..,b_{i-1},b_{i+1},..,b_{n})$ for each $i$. This is equivalent to $b_{i}\notin acl(A,b_{1},..,b_{i-1})$ for all $i=1,..,n$.  For $\bar b$ a tuple from $D$, $dim({\bar b}/A)$ denotes the cardinality of some/any maximal algebraically independent over $A$ subtuple of $\bar b$. And it turns out that for a subset  $X$ of $D^{n}$, definable over $A$, $dim(X)$ as in Definition 1.4 (i) coincides with $max\{dim({\bar b}/A):{\bar b}\in X\}$. 
There is  a unique type over $A$ of any algebraically independent $n$-tuple $(b_{1},..,b_{n})$. And this explains Fact 1.5(ii) for example.  

What we have called dimension is the same thing as Morley rank (or $RM$) and what we called multiplicity is the same thing as Morley degree $dM$, all in the special case of strongly minimal theories.  So we will use these expressions freely.  We also require: 

\begin{Fact}
 Suppose $X\subseteq D^{n}$ is definable, and $Y\subseteq D^{m}$ is definable. Then $RM(X\times Y) = RM(X) + RM(Y)$ and $dM(X\times Y) = dM(X).dM(Y)$.  In particular when $X$ and $Y$ both have Morley degree $1$ so does $X\times Y$. 
\end{Fact} 

\vspace{5mm}
\noindent
We now say a few words about our formalism for pseudofinite theories and structures.  The usual way of considering pseudofinite structures is via ultraproducts: up to elementary equivalence a pseudofinite structure is an ultraproduct of finite structures, and as such lives in a nonstandard model of set theory. See section 3 of \cite{Tao} for a nice account. Another way of doing this, avoiding explicit reference to ultraproducts is as follows: Let $M$ be a ($1$-sorted)  pseudofinite structure which is {\em saturated} ($\kappa$-saturated and of cardinality $\kappa$ for some large $\kappa$). Let $T$ be the theory of $M$, in its language $L$.  For each $L$-sentence $\sigma\in T$, let  $M_{\sigma}$ be a finite model of $\sigma$.  Let us fix such $M_{\sigma}$ and expand $M_{\sigma}$ to an $L'$-structure $M_{\sigma}'$ for a suitable language $L'$ containing $L$ which I will describe now. $L'$ contains additional sorts for the power set ${\cal P}(M_{\sigma}^{n})$ of each Cartesian power of $M_{\sigma}$, as well as the appropriate membership relation. $L'$ also has sorts for the natural numbers $\N$ and real numbers $\R$ and for the usual embedding of $\N$ into $\R$ as well as all arithmetic operations.  We also have for each $n$ a symbol for a function from ${\cal P}(M_{\sigma}^{n})$ to $\N$, the intended interpretation being cardinality. Let $M_{\sigma}'$ be the tautological expansion of $M_{\sigma}$ to an $L'$-structure. We can of course do this uniformly for each $\sigma\in T$ to get a family $\cal K = \{M_{\sigma}':\sigma\in T\}$ of $L'$-structures.  We note that $T \cup Th({\cal K})$ is consistent, hence has a saturated model $M'$ say of cardinality $\kappa$. The $L$-part or $L$-reduct of $M'$ is a saturated model of $T$ so can be identified with $M$. Hence $M'$ is an expansion of $M$ to a model of $Th(\cal K)$.  Note that among the subsets of $M^{n}$ we have (i) the sets which are definable in the language $L$ with parameters from $M$, and (ii) more generally the sets which are definable in the structure $M'$ with parameters, which are precisely the sets corresponding to the the interpretation in $M'$ of the sort for sets of $n$-tuples. We call sets of kind (ii) {\em internal}, and any internal set $X$ has  by virtue of the cardinality functions in $L'$ a nonstandard finite cardinality which we write as $|X|\in {\N}^{*}$  (the interpretation of the sort for $\N$ in the structure $M'$).  Of course sets of kind (i) are among the internal sets but have a more privileged status from the point of view of the model theory of the first order theory $T$. 

Note that there are some choices here: the choice of a completion of $T\cup Th({\cal K})$ is analogous to a choice of an ultrafilter in the usual presentation. 

\section{Proofs and additional remarks} 

%Recall the basic notions. A structure $M$ is in language $L$ is pseudofinite if every sentence true in $M$ is true in some finite $L$-structure.  Equivalently $M$ is elementarily equivalent %to an ultraproduct of finite $L$-structures. If $M$ is pseudofinite and saturated say, then every definable set $X$ in $M$ has a ``nonstandard finite cardinality"  $|X|$ which is an $element of a saturated elementary extension of $(\N,+,\times,<,....)$, and the map taking $X$ to $|X|$ satisfies the usual properties inherited from the finite setting. 

%Suppose $D = M$ is strongly minimal and saturated. $D$ is said to be unimodular if whenever $a = (a_{1},..,a_{n})$ and $b = (b_{1},..,b_{n})$ are each independent $n$-tuiples from %$D$ and $a\in acl(b)$  (so also $b\in acl(a)$)  then $mlt(a/b) = mlt(b/a)$. 

%Definable means possibly with parameters. We refer to \cite{Pillay-book} for basics of stability, Morley rank ($RM(-)$) etc. 

{\em Proof of Theorem 1.1.}
\newline
(i) We prove it by induction on $RM(X)$.  When $RM(X) = 0$, $X$ is finite so $|X| = |X|$. 
Now suppose $RM(X) = n\geq 1$  where  $X$ is a definable subset of $D^{m}$ (for some $m\geq n$).  By Fact 1.5 (iii), we may assume that $X$ is an $n$-cell.  So 
under some projection from $D^{m}$ to $D^{n}$, $\pi(X)$ has Morley  rank $n$ and $\pi|X$ is $r$ to one for some positive (standard) natural number $r$.    So $|X| = t|\pi(X)|$. 
And $|\pi(X)|  = |D^{n}| - |D^{n}\setminus \pi(X)|$. Now $|D^{n}| = q^{n}$, and $RM(D^{n})\setminus \pi(X)$ has Morley rank $< n$. So we can apply the induction hypothesis to see that $|\pi(X)|  = q^{n} - R(q)$ where $R$ is a polynomial over $\Z$ with degree $< n$.  So $|X| = rq^{n} - rR(q)$ and we are finished. 

\vspace{2mm}
\noindent
(ii)   We start with
\newline
{\em Claim.} For any $L$-formula $\phi({\bar x}, {\bar y})$ where ${\bar x} = (x_{1},..,x_{n})$, and ${\bar a}$ (for the ${\bar y}$ variables), there is an $L$- formula $\psi({\bar y}) \in tp({\bar a})$ such that for all ${\bar b}$ satisfying $\psi$, $P_{\phi({\bar x},{\bar a})}(q) =  P_{\phi({\bar x},{\bar b})}(q)$. 
\newline
{\em Proof of claim.} 
Write $\phi({\bar x},{\bar a})$ as a ``disjoint union" of (formulas defining) cells $\phi_{1}({\bar x},{\bar a_{1}})$,...., $\phi_{s}({\bar x},{\bar a}_{s})$. 
Fix one of these formulas, say $\phi_{1}({\bar x},{\bar a_{1}})$.  So the solution set of this formula projects $t$-to-one (some $t$) to coordinate axes  ${\bar x}'$ and the projection $\psi({\bar x}',{\bar a_{1}})$  is a definable subset of ${\bar x}'$-space of maximal Morley rank. Apply induction to the formula $\neg\psi({\bar x}',{\bar a_{1}})$.  Likewise for each of the other $\phi_{i}({\bar x},{\bar a_{i}'})$.  Details are left to the reader, bearing in mind the proof of part (i).

\vspace{2mm}
\noindent
We can apply compactness to the claim to obtain the desired conclusion. 

\vspace{5mm}
\noindent
{\em Proof of Corollary 1.2.}
\newline
Here we start to write  tuples from $D$ as $a,b,etc..$ hopefully with no confusion. 
  Let $a,b\in D^{n}$ each be an $n$-tuple of algebraically independent (over $\emptyset$) elements of $D$, such that $a$ and $b$ are interalgebraic.
 Let $k = mult(b/a)$ and $\ell= mult(a/b)$ as defined in the paragraph following the statement of Theorem 1.1 in section 1. We have to prove that $k=\ell$. 
 Let $\psi(x,y)$ be an $L$-formula such that $D \models \phi(a,b))$,  $\psi(a,y)$ isolates $tp(b/a)$ and $\psi(x,b)$ isolates $tp(a/b)$. Let $\phi_{1}(x)$ be $\exists ^{=k}y(\psi(x,y))$ and $\phi_{2}(y)$ be $\exists^{=\ell} x(\psi(x,y)$. Let $\chi(x,y)$ be the formula $ \phi(x,y)\wedge\phi_{1}(x)\wedge\phi_{2}(y)$. So $\chi(x,y)$  is true of $(a,b)$ in $D$.  Let $Z\subseteq D^{2n}$ be the set defined by $\chi(x,y)$. We compute $|Z|$ in two ways.  Let $X$ be the projection of $Z$ on the first $n$-coordinates, and $Y$ the projection of $Z$ on the last $n$ coordinates.
Then $|Z| = k|X| = \ell|Y|$.  Note that $X$ and $Y$ have Morley rank $n$ hence by Fact 1.5 and Theorem 1.1, there are polynomials $P(x)$, $Q(x)$ over $\Z$ of degree $< n$ such that $|X| = q^{n}-P(q)$ and 
$|Y| = q^{n} - Q(q)$.  If by way of contradiction $k>\ell$ we have $(k-\ell)(q^{n})  = kP(q)-\ell Q(q)$.   This is impossible, as the right hand side is an integral polynomial of degree $< n$ in the infinite nonstandard natural number $q$.
\newline
Thus proves unimodularity of $D$. Local modularity follows as remarked earlier, but we will discuss further below local modularity (including the definition).

\vspace{5mm}
\noindent
{\em Proof of Corollary 1.3.}
\newline
By the assumptions and Fact 1.6, $RM(V\times W) = d_{1} + d_{2}$, hence $P_{V\times W}$ has degree $d_{1} + d_{2}$ by Theorem 1.1, and also   $V\times W$ has Morley degree $1$.  So either (i) $RM(E) = d_{1} + d_{2}$ in which case $RM((V\times W)\setminus E) < d_{1} + d_{2}$, or (ii) $RM(E) < d_{1} + d_{2}$.  In case (i), $|(V\times W)  \setminus E| = R(q)$ for $R$ an integral polynomial with degree $< d_{1} + d_{2}$ whereby 
$|E| = |V||W|-R(q)$, i.e. all but $R(q)$ elements of $V\times W$ are in the edge relation $E$.  It follows immediately that for any subsets $A$ of $V$ and $B$ of $W$ which are internal (so $A\times B$ has a well-defined nonstandard finite cardinality), all but at most $R(q)$ elements of $A\times B$ are in the edge relation $E$.  Similarly in case (ii).  End of proof.

\vspace{5mm}
\noindent
I will now discuss some other consequences and questions. At this point I will assume more familliarity with stability theory, using say \cite{Pillay-book} as a reference. 
Recall that local modularity of a strongly minimal set $D$ means that, after naming some parameters, for finite tuples $b, c$ from $D$, 
$dim(b,c) = dim(b) + dim(c) + dim(acl(b)\cap acl(c))$. It is equivalent to  $D$ being $1$-based: namely for all tuples $b,c$ in $D$ $b$ and $c$ are independent over $acl^{eq}(b) \cap acl^{eq}(c)$ where independence is in the sense of nonforking.  But $1$-basedness now makes sense for any stable theory. 

Now suppose that $T$ is a complete theory of finite $U$-rank, in the sense that every complete type has finite $U$-rank. Suppose also that $T$ is pseudofinite. Then by Corollary 1.2, any strongly minimal definable set in $T$ is locally modular.  As in the proof of Proposition 3.5 of \cite{KP} we conclude:
\begin{Corollary} Any pseudofinite theory of finite $U$-rank is $1$-based. In particular any group definable in such a theory is abelian-by-finite. 
\end{Corollary} 

In \cite{MT} it is proved that any stable pseudofinite group $G$ is solvable-by-finite, and some examples of such $G$ which are not even nilpotent-by-finite are mentioned.  We see in particular that a stable pseudofinite theory need not be $1$-based. It is reasonable to ask what can be said in general about stable pseudofinite theories.

\begin{Problem} Show that in a stable pseudofinite theory, every regular type is locally modular. 
\end{Problem}

It may also  be of interest to give a conceptual pseudofinite account of the Malliaris-Shelah stable regularity lemma \cite{MS}.  The model-theoretic context is not a pseudofinite stable theory, but rather
 a pseudofinite bipartite graph $(V,W, E)$, living in a nonstandard model of a set theory,  where the  formula $xEy$ stable.  

\vspace{2mm}
\noindent
It is also natural to ask what is the appropriate level of generality of the precise counting result in Theorem 1.1.   Firstly there should be no problem obtaining a similar result for pseudofinite $\aleph_{1}$-categorical theories, where again any definable set will have cardinality an integral polynomial in $q$ where $q$ is the cardinality of a given strongly minimal set.  Likewise for pseudofinite  theories of finite $U$-rank ``coordinatized" by finitely many strongly minimal sets  $D_{1},..,D_{r}$ where the  relevant polynomials would be in $(q_{1},..,q_{r})$.
It would possibly be interesting (although laborious) to consider appropriate statements and proofs in the case of pseudofinite theories of finite $U$-rank, with Morley rank equal to $U$-rank, and Morley rank definable.

Finally let us remark that the classical results of Zilber and Cherlin, Harrington, Lachlan (see Chapters 2 and 3 of \cite{Pillay-book})  say that strongly minimal (in fact $\omega$-stable) $\omega$-categorical theories are pseudofinite.

\end{document}